\documentclass[a4paper,oneside,english,reqno]{amsart}

\usepackage[utf8]{inputenx}
\usepackage{amsthm, amsmath, amssymb, amsfonts, mathrsfs, mathtools, stmaryrd}
\usepackage[english]{babel} 
\usepackage{textcomp, url, enumerate, latexsym, graphicx, varioref, hyperref, multirow, layout, siunitx, booktabs}
\usepackage{pdflscape}
\usepackage{verbatim}
\usepackage{geometry}
\usepackage{dsfont}
\usepackage{csquotes}

\usepackage[pagewise,displaymath,mathlines]{lineno}

\RequirePackage{enumitem}
\setlist[itemize]{font = \upshape, before = \leavevmode}
\setlist[enumerate]{font = \upshape, before = \leavevmode}
\setlist[description]{before = \leavevmode}

\usepackage[usenames]{color}
\usepackage{colortbl}

\usepackage{tikz, babel, rotating, tikz-cd, pigpen}
\usetikzlibrary{matrix, arrows, decorations.pathmorphing}
\usepackage{cleveref}
\RequirePackage[color       = blue!20,
                bordercolor = black,
                textsize    = tiny,
                figwidth    = 0.99\linewidth]{todonotes}
\usepackage[all]{xy}
\usepackage{microtype}

\RequirePackage[backend    = biber,
                sortcites  = true,
                giveninits = true,
                doi        = false,
                isbn       = false,
                url        = false,
                maxnames   = 6,
                style      = numeric]{biblatex}
\DeclareFieldFormat*{title}{#1}
\renewbibmacro{in:}{}

\numberwithin{equation}{section}

\numberwithin{equation}{section}

\theoremstyle{plain}
\newtheorem{theorem}[equation]{Theorem}

\newtheorem{corollary}[equation]{Corollary}
\newtheorem{proposition}[equation]{Proposition}
\newtheorem{lemma}[equation]{Lemma}
\crefname{lemma}{Lemma}{Lemmas}

\newtheorem{thmX}{Theorem}

\newtheorem{factX}{Fact}

\theoremstyle{definition}

\newtheorem{definition}[equation]{Definition}

\theoremstyle{remark}
\newtheorem{remark}[equation]{Remark}

\newcommand{\A}{\mathbb{A}}
\newcommand{\PP}{\mathbb{P}}
\newcommand{\DM}{\mathbf{DM}}
\newcommand{\Sm}{\mathrm{Sm}}
\newcommand{\EssSm}{\mathrm{EssSm}}
\newcommand{\nis}{\mathrm{nis}}
\newcommand{\Gm}{\mathbb{G}_m}

\newcommand{\zar}{\mathrm{Zar}}
\newcommand{\Pic}{\mathrm{Pic}}
\newcommand{\Cor}{\mathrm{Cor}}

\newcommand{\Iso}{\mathrm{Iso}}
\newcommand{\Sch}{\mathrm{Sch}}
\newcommand{\perf}{\mathrm{perf}}

\newcommand{\Spec}{\operatorname{Spec}}
\newcommand{\Fieldsop}{\mathrm{Points}}

\newcommand{\bbZ}{\mathbb Z}
\newcommand{\codim}{\operatorname{codim}}

\newcommand{\pt}{\mathrm{pt}}
\newcommand{\Hom}{\mathrm{Hom}}

\makeatletter
\@namedef{subjclassname@2020}{%
  \textup{2020} Mathematics Subject Classification}
\makeatother

\author{Druzhinin A.E., Urazbaev A.A.}
\title{
Non-equivalences of motivic 
codimension filtration quotients.
}

\address{Andrei Druzhinin, \\
Chebyshev Laboratory, St. Petersburg State University, 14th Line V.O., 29, Saint Petersburg 199178 Russia
}

\address{Askar Urazbaev, \\
St. Petersburg State University, Department of Mathematics and Computer Sciences, 14th Line V.O., 29, Saint Petersburg 199178 Russia
}


\thanks{Research is supported by the Russian Science Foundation grant 19-71-30002}

\begin{document}


    
    

\begin{abstract}
    We prove that 
    a
    motivic equivalence 
    of objects of the form 
    \begin{equation*}
    X/(X-x)\simeq X^\prime/(X^\prime-x^\prime)
    \end{equation*}
    in 
    $\mathbf{H}^\bullet(B)$
    or
    $\mathbf{DM}(B)$
    over a scheme $B$,
    where 
    $x$ and $x^\prime$ are closed points of smooth $B$-schemes 
    $X$ and $X^\prime$,
    implies
    an isomorphism of residue fields, i.e.
    \[x\cong x^\prime.\]
    
    For a given $d\geq 0$,
    $X,X^\prime\in\mathrm{Sm}_B$, 
    $\operatorname{dim}_B X=d=\operatorname{dim}_B X^\prime$,
    and 
    closed points 
    $x$ and $x^\prime$ 
    that residue fields are simple extensions 
    of the ones of $B$,
    we show an
    isomorphism of groups
    \[\Hom_{\mathbf{DM}(B)}(X/(X-x),X^\prime/(X^\prime-x^\prime)))\cong\Cor(x,x^\prime),\]
    and prove that it leads to an equivalence of subcategories.
    
    Additionally, 
    using
    the result on perverse homotopy heart 
    by F.~Déglise and N.~Feld and F.~Jin
    and 
    the strict homotopy invariance theorem for presheaves with transfers
    over fields
    by the first author,
    we
    prove
    an equivalence
    of the Rost cycle modules category
    and the homotopy heart of $\mathbf{DM}(k)$
     over a field $k$ with integral coefficients.

\end{abstract}

\maketitle

\section{Introduction}
We study 
some analogues 
of topological spheres 
in 
the pointed Morel-Voevodsky motivic homotopy 
category $\mathbf{H}^\bullet(B)$ \cite{MV,Voe98}, the stable one 
$\mathbf{SH}(B)$ \cite{Jardine-spt,Voe98}, \cite[Appendix C]{Hoyois_SHbaseschemeBchN}, and the Voevodsky motives category $\DM(B)$  \cite{Voe-motives,Cisinski-Deglise-Triangmixedmotives}
over a base scheme $B$.
%
The "only if" parts of
our classifiaction results,
Theorems \ref{thmX:th-postintro:classif:cls} and \ref{th-postintro:classif},
are new for any base field or base scheme,
while 
the "if" parts
are 
well-known 
provided by 
\Cref{fact:ClassifParta}.
%
Theorems 
\ref{thmX:eq:hearttshomottructure} and \ref{th-intro:dist}
are 
known 
over base schemes with perfect residue fields or after inverting of exponential characteristics.

Topological
spheres 
$S^d$ 
generate 
the category $\mathbf{H}^\mathbf{top}$
via extensions,
and
the family of functors
\begin{equation}\label{eq:HomSdHtopSetbullet}
\pi_{d}(-)=\Hom_{\mathbf{H}^\mathbf{top}}(S^d,-)\colon \mathbf{H}^\mathbf{top} \to \mathrm{Set}^\bullet, 
\quad d\in \mathbb Z_{\geq 0},
\quad\text{see Notation 
\ref{notation:pointedcategories},}
\end{equation}
is conservative,
because for 
a CW-complex $X$, 
there is a filtration by open 
subspaces 
\begin{gather}
\label{eq:filtrtop}
X \supset \dots \supset X^{c} \supset X^{c-1}\supset \dots \supset X^{0}
\\
\label{eq:XdXd-1simeqcoprodSd}
X^{c}/X^{c-1}\simeq \bigvee_{C_d} S^d, \quad C_c\in\mathrm{Set}, \; c\in\mathbb Z_{\geq 0}.
\end{gather}
For any scheme $B$,
the category $\mathbf{H}^\bullet(B)$
is not generated
by $S^d$ and $\Gm^{\wedge d}$
via smash products and extensions, 
see Notation \ref{notation:SdGmw1} and \ref{notation:pointedcategories}.
For $X\in\Sm_B$ in general,
there is no filtration 
like \eqref{eq:filtrtop}  
satisfying condition \eqref{eq:XdXd-1simeqcoprodSd} 
with $S^d$ substituted by
motivic spheres
\begin{equation}\label{eq:sphereTd}T^{\wedge d} = 
= S^d\wedge\Gm^{\wedge d}\simeq\mathbb A^d_B/(\mathbb A^d_B-0_B), 
\quad d\in\mathbb Z_{\geq 0}.\end{equation}
An analogue
called \emph{codimension filtration}
is the limit of filtrations 
\eqref{eq:filtrtop}
such that 
$\codim_X Z^{c}=c$,
$Z^c=X\setminus X^{c}$. 
While
$S^d$ are ``building blocks'' of a topological space $X$, 
such ``building blocks'' of $X\in\Sm_B$ are 
motivic pro-spaces of the form 
\begin{equation}\label{eq:Xx/Xx-x}
    X_x/(X_x-x),\quad\text{see Notation \ref{notation:propobjcofib},}
\end{equation}
where $x\in X$ is a point, and $X_x=\Spec \mathcal O_{X,x}$ is the local scheme of $X$ at $x$.
Pro-spaces \eqref{eq:Xx/Xx-x} 
control 
motivic homotopy types,
in the sense that
the family of functors 
\begin{equation}\label{eq:intro:conservativefamily:HBXXx}
\Hom_{\mathbf{H}^\bullet(B)}(X_x/(X_x-x), -)\colon \mathbf{H}^\bullet(B)\to \mathrm{Set}^\bullet,
\quad X\in\Sm_B, x\in X,
\end{equation}
is conservative.
Note that 
a
pro-object \eqref{eq:Xx/Xx-x} 
is an object 
for $x\in X$ being a closed point.



\subsection{Classification}

The Grothendieck six-functors formalism 
provided by \cite{zbMATH05292568,zbMATH05318528,Cisinski-Deglise-Triangmixedmotives}
implies the following result.

\begin{factX}\label{fact:ClassifParta}
Let 
$p\colon X\to B$, $p^\prime\colon X^\prime\to B$ be essentially smooth morphisms of schemes,
$X$ and $X^\prime$ be local, 
${x}$ and ${x}^\prime$ be closed points.
Suppose that 
\begin{equation}\label{eq:dimsimpareequal:facta}
\dim_{B} X = \dim_B X^\prime 
\in \mathbb Z,\quad  
{x}\cong {x}^\prime\in \Sch_B,
\end{equation}
then there is an isomorphism 
of pro-objects
in 
$\mathbf{SH}(B)$ 
\begin{equation}\label{eq:XxXxxsimeqprimeone}
X/(X-{x})\simeq X^\prime/(X^\prime-{x}^\prime).
\end{equation}
\end{factX}

\begin{remark}
    \Cref{fact:ClassifParta}
    follows
    because of isomorphisms
    $X/(X-\mathrm{x})\simeq \mathbf{x}_!({x})$, $\mathbf{x}\colon {x}\to X$.
\end{remark}

\begin{remark}\label{rem:Nisnevicheq}
    The recent 
    preprint \cite{NisEqSmOP_DrUr}
    presents
    a short elementary proof of the equivalence
    \eqref{eq:XxXxxsimeqprimeone}
    in the category $\mathbf{H}^\bullet(B)$,
    and moreover, the category of pointed Nisnevich sheaves.
\end{remark}


To complete the classification of \eqref{eq:Xx/Xx-x} means to prove 
the converse implication in \Cref{fact:ClassifParta}.
%
The main result of this paper is as follows.

\begin{thmX}[\Cref{th-postintro:classif:cls}]\label{thmX:th-postintro:classif:cls}

Let $X,X^\prime\in\Sm_B$,
and let $x\in X$, $x^\prime\in X^\prime$ be closed points.
Then an isomorphism of objects
\begin{equation*}
X/(X-\mathrm{x})\simeq X^\prime/(X^\prime-\mathrm{x}^\prime)
\end{equation*}
in 
$\mathbf{H}^\bullet(B)$, 
$\mathbf{SH}(B)$ 
or $\mathbf{DM}(B)$
holds if and only if 
\begin{equation}\label{eq:into:claim:xsimeqxprime}
\mathrm{x}\cong \mathrm{x}^\prime\in \Sch_B,
\end{equation}
and
$\dim^{x}_{B}X=\dim^{x^\prime}_B X^\prime$.
\end{thmX}



\Cref{thmX:th-postintro:classif:cls} completes the classification 
for $x$ being closed points.
Also we have a partial result
for pro-objects \eqref{eq:Xx/Xx-x}
with non-closed points $x$ in $X\in\Sm_B$,
see \Cref{th-postintro:classif}.

\subsection{Additional results on cycle modules and $\mathrm{Hom}$-groups} 

The mapping
\begin{equation*}
    E\longmapsto \big((l,x) \mapsto \mathrm{Hom}(X/(X-x),E\{l\})\big),
\end{equation*}
see Notation \ref{notation:SdGmw1},
gives rise to the functor
\begin{equation}\label{eq:intro:DMtoMWmod}
\mathbf{DM}(k)
\to
\mathcal{MC}ycle(k),
\end{equation}
where
the right side is
the category of
Rost
cycle modules
introduced in
\cite{Rost1996}.
It is proven in \cite{DegliseModuleshomotopiques,déglise2022perversehomotopyheartmwmodules}
that the functor \eqref{eq:intro:DMtoMWmod}
induces an equivalence on homotopy hearts
for a perfect base field $k$ or after inverting of the exponential characteristic
in the categories.
In \Cref{sect:CycleModules},
we show that combynination \cite{DegliseModuleshomotopiques,déglise2022perversehomotopyheartmwmodules} and \cite{StrictA1invariancefields}
allows to
extend the generality.
\begin{thmX}[\Cref{eq:hearttshomottructure}]\label{thmX:eq:hearttshomottructure}
    For any field $k$,
    there is the equivalence
    \begin{equation}\label{eq:MmodDMheart:intro}
    \mathcal{MC}ycle(k)
    \simeq
    \mathbf{DM}(k)^{\heartsuit}
    \end{equation}
    where $\heartsuit$ relates to the homotopy t-structure heart.
\end{thmX}
\begin{remark}
In
the context
of equivalence
\eqref{eq:MmodDMheart:intro},
\Cref{thmX:th-postintro:classif:cls} and 
\Cref{th-postintro:classif} 
intend to show
that
the category 
$\mathcal{MC}ycle(k)$
can not be reduced to 
a smaller one.
\end{remark}

Whenever we have 
equivalence \eqref{eq:MmodDMheart:intro}
it is natural to ask 
what is the cyclemodule that corresponds to
the object \eqref{eq:Xx/Xx-x}
along 
\eqref{eq:intro:DMtoMWmod}.
In other terms
this means 
to compute 
Hom-groups
in $\DM(k)$
between objects
and pro-objects
\eqref{eq:Xx/Xx-x} and its suspensions.
We have the following result in this direction.

\begin{thmX}[\Cref{th:onedim:subcatDMCorB}]\label{th-intro:dist} 

Let $B$ be a scheme as in Notation \ref{notation:categoriesSH(B)DM(B)}.

For each $d\in \mathbb Z_{\geq 0}$,
and $(X,x),(X^\prime,x^\prime)\in\Sm^{d,\cdot,1}_B$,
see Notation  \ref{notation:Smddot},
there is an isomorphism of groups
\[\Hom_{\DM(B)}(X/(X-x),X^\prime/(X^\prime-x^\prime)[l]))\cong\begin{cases}
    \Cor(x,x^\prime),&l=0,\\
    0,&l\neq 0,
\end{cases}\]
and
moreover,
there is 
an equivalence of 
the subcategories in $\DM(B)$ and $\mathbb Z\times \Cor(\Sch_B)$
\begin{equation*}
\begin{array}{lcl}
(\;
X/(X-x)[l]\;|\;
    (X,x)\in\Sm^{d,\cdot,1}_B, l\in \mathbb Z
\;)_{\DM(B)}
&\simeq&
\mathbb Z\times \Cor_B(\Fieldsop^1_B),\\
X/(X-x)[l]&\mapsto& (l,x)
\end{array}
\end{equation*}
where $\mathbb Z$ denotes the discrete category that objects are integers, see Notation 
\ref{notation:PointsB} and \ref{notation:fsubcat}.

\end{thmX}

\subsection{Notation and conventions}
\begin{enumerate}

\item \label{notation:SdGmw1}
$S^d=\Delta^d/\partial \Delta^d$,
$\Gm^{\wedge 1}=\Gm/\{1\}$.
$E[d]=E\wedge S^d$,
$E\{d\}=E\wedge\Gm^{\wedge d}$.
\item \label{notation:Twd}
$T^{\wedge d} = S^d\wedge\Gm^{\wedge d}$ that is motivically equivalent to
$\mathbb A^d_B/(\mathbb A^d_B-0_B)$.

\item
$\Sch_B$, $\Sm_B$, $\EssSm_B$  are the category of $B$-schemes, smooth and essentially smooth schemes,
which are limits of smooth ones with \'etale affine transition maps.

\item
$X-Z$ and $X\setminus U$
are the open and the reduced closed complements
for a closed $Z$ and an open $U$ in $X$ respectively.

\item \label{notation:pointedcategories}
$\mathrm{Set}$, $\mathbf{H}^\mathrm{top}$, $\mathbf{H}(B)$ 
are
the category of sets, 
the homotopy category of topological spaces, 
the Morel-Voevodsky motivic homotopy category over a scheme $B$.
Let
$\mathrm{Set}^\bullet$, $\mathbf{H}^\mathbf{top}$, 
$\mathbf{H}^\bullet(B)$ are the respective pointed ones. 

\item 
$(-)_+\colon X\mapsto X_+:=X\amalg *$ is the functor from an unpointed category to the pointed one.

\item \label{notation:categoriesSH(B)DM(B)}
We write $\mathbf{SH}(B)$ 
for the Voevodsky stable motivic category 
in sense of
\cite{MV,Voe98,Jardine-spt}
or
\cite[Appendix C]{Hoyois_SHbaseschemeBchN}.
Similarly 
for 
$\mathbf{DM}(B)$ 
whenever it is well defined,
see
\cite{Voe-motives,Cisinski-Deglise-Triangmixedmotives}.

\item
Given a morphism 
$f\colon S\to B$,
$f^*$ and $f_*$
denote the inverse and direct image functors.


\item \label{notation:hocofibdoublequotientsymbol} 
Given 
$Y\to X$,
symbols
$X/Y$, 
$\operatorname{hocofib}(Y\to X)$ 
denote 
the cofiber, 
the homotopy cofiber.






\item \label{notation:propobjcofib}
Given pro-objects 
$X=\varprojlim_{\alpha\in A} X_\alpha$ and $U=\varprojlim_{\alpha\in A} U_\alpha$ 
for a filtered set $A$, 
and a morphism 
$U\to X$
defined by morphisms $U_\alpha\to X_\alpha$, 
\[X/U=\varprojlim_{\alpha\in A}X_\alpha/U_\alpha.\]


\item
Given a smooth scheme over $B$,
the same symbol denotes 
the representable presheaf.

\item
Given a pro-presheaf over $B$,
the same symbol denotes
the 
pro-object in
$\mathbf{H}(B)$ or $\DM(B)$. 


\item 
We write  $Z\not\hookrightarrow X$ for a closed immersion of schemes.


\item 
—, and a point $x\in X$, we denote by $X_x$ the local scheme $\Spec \mathcal O_{X,x}$.

\item\label{notation:Xpcp}
$X^{(c)}=\{x\in X|\codim_X x=c\}\in\mathrm{Set}$, 
and $X^{(c)}=\coprod_{x\in X^{(c)}}x\in\Sch$.
\item\label{notation:Xpdp}
$X_{(d)}=\{x\in X|\codim_X x=\dim^x X-d\}\in\mathrm{Set}$, 
and $X_{(d)}=\coprod_{x\in X_{(d)}}x\in\Sch$.

\item
Given a scheme $V\in \Sm_B$, 
use notation $V/\pt_B=\mathrm{Cone}(V\to \pt_B)[-1]\in \DM(B)$.

\item 
Denote by $L_\nis$ the endofunctor on the derived category of abelian presheaves with transfers
given by the Nisnevich local replacement, and write $L_\zar$ for the Zariski one.

\item 
We denote hom-groups in $\DM(k)$ by
$\Hom_{\DM(k)}(-,-)$ or $\Hom(-,-)$.

\item\label{notation:hmotlsimeqHomDMk}
We use notation 
$h^l_\mathrm{mot}(F)$
for the presheaf on $\Sm_k$
given by
$\Hom_{\mathbf{DM}(k)}(-,F[l])$.

\item 
$Z_d(X)$ denotes the group of $d$-dimensional cycles in a scheme $X$.

\item \label{notation:divf}
$\operatorname{div}(f)$ denotes the divisor of an invertible rational function $f$.
\item \label{notation:divelementinCor}
Given $U,V\in\Sch_B$, $\dim_B V=1$, 
write $\operatorname{div}(f)\in\Cor(V_0,V_1)$ 
for
the respective element.
\item \label{notation:cCortoHom}
Given $c\in\Cor_B(U,V)$, denote by the same symbol the image in $\Hom_{\DM(k)}(U,V)$.

\item\label{notation:zperfzalg}
Given $z=\Spec k$ for a field $k$, 
denote
$z^\perf = \Spec k^\perf$,
$z^\mathrm{alg} = \Spec k^\mathrm{alg}$,
where
$k^\perf$ and $k^\mathrm{alg}$
are the purely inseparable closure and the algebraic closure of $k$.

\item\label{notation:sdeg}
%
$\operatorname{sdeg}_{z}x= \operatorname{sdeg}_{K}L$
is the separable degree 
for $x=\Spec L$ and $z=\Spec K$.

\item 
$\Iso(C)$ or $\Iso C$
is the set of the isomorphism classes of the objects of 
a small category $C$.

\item \label{notation:fsubcat}
$F_{C}=(E\in F)_{C}=(E| E\in F )_{C}$
is the subcategory of $C$ spanned a family of objects $F$.



\item 
$\Cor(C)\subset\Cor(\Sch_B)$ 
is the subcategory
spanned by 
the objects of $C$.

\item\label{notation:FunctorsCat}
Given a category $C$ and a small category $D$, 
$C^D$ denotes the category of functors $D\to C$. 
$\Delta_1$ denotes the full subcategory spanned by $[0]$ and $[1]$ in $\Delta$.

\item\label{notation:Corpair}
$\Cor_B^\mathrm{pair}$ and 
$\Cor_B^\mathrm{c.o.pair}$
are the full subcategoryies in $\Cor_B^{\Delta_1^\mathrm{op}}$
spanned by open immersions and clopen immersions in $\Sm_k$. 
Given $X\in\Sm_k$, and open subscheme $U$ in $X$, we write $(X,U)=(U\hookrightarrow X)\in\Cor_B^\mathrm{pair}$.  
Similarly for $\Cor^\mathrm{pair}(\Sch_B)$ and 
$\Cor^\mathrm{c.o.pair}(\Sch_B)$
\item\label{notation:CorCorpair(P0P1)}
For any objects $P_0,P_1\in\Cor_B^\mathrm{pair}$, we write 
$\Cor_B(P_0,P_1)=\Cor_B^\mathrm{pair}(P_0,P_1)$.

\item
$S^\amalg=(S)^\amalg\subset C$ is 
the coproduct closure of
a subcategory $S\subset C$. 

\item \label{notation:PointsB}
$\Fieldsop_B$
is the category of irreducible reduced zero-dimensional finite $B$-schemes,
and
$\Fieldsop^l_B$ 
is 
spanned by the ones 
that residue fields extensions has $l$ generators. 


\item\label{notation:Smddot} 
$\Sm^\cdot_B$ is the category of pairs $(X,x)$, 
where $X\in \Sm_B$, $x\in X$ is a closed point;\\
$\Sm^{d,\cdot}_B$ and $\Sm^{\cdot,l}_B$ 
are the subcategories of such ones that 
$\dim^x_B X=d$
and
$x\in \Fieldsop^l_B$
respectively;\\ 
$\Sm^{d,\cdot,l}_B=\Sm^{d,\cdot}_B\cap\Sm^{\cdot,l}_B$.

\item\label{notation:catZ}
$\mathbb Z$ denotes the discrete category that objects are integers.




\end{enumerate}

\section{One-dimensional open pairs}\label{sect:MotCohGenCircle}

In this section, 
we study 
motives 
\[\begin{array}{lllllll}
T_{k(z)} &=& z^*(T)&=&\A^1_{k(z)}/(\A^1_{k(z)}-0)&\in&  \mathbf{DM}(k(z)),\\
\widetilde T^z &=& &&\A^1_k/(\A^1_k-z)&\in& \mathbf{DM}(k),
\end{array}\]
for closed points $z\in \A^1_k$,
where 
$z^*$ and $z_*$ stand for 
the inverse and direct images along the base change functor $\Sm_k\to\Sm_z:X\mapsto X\times_k z$.
\begin{remark}
For $k(z)/k$ being purely inseparable,
\begin{equation}\label{eq:zuswTzsimeqzusT}z^*(\widetilde{T}^z)\simeq T_{k(z)},\end{equation}
because $(z\times_k z)_\mathrm{red}\simeq z$.
\end{remark}
The results:
\begin{itemize}
    \item 
    \Cref{cor:extensioncircle}
    generalises computation of the weight one motivic homologies $H^{p,1}(-,\bbZ)$ \cite[Corollary 4.1]{MVW};
    \item 
    \Cref{prop:CatC/El} dscribes the homomorphims 
    in between of the objects $\widetilde T^z$;
    \item
    \Cref{cor:prop:CatC/E} 
    shows some kind of ``weak $\bbZ/p\bbZ$-orthogonality'' 
    of $T$ and $\widetilde T^z$ 
    for purely inseparable extensions $k(z)/k$,
    where $p=\operatorname{char} k$. 
\end{itemize}
\begin{remark}
    Though this is not used, let us note that, 
    if $k$ is perfect the results are trivial or follows easily from the results of \cite{MVW}.
\end{remark}

\subsection{Homomorphisms from a scheme to a pair}


\begin{theorem}\label{lm:LA1htr(V)}
Let $V = \A^1_k-D$ for some closed subscheme $D$ in $\A^1_k$ of positive codimension.
There is a quasi-isomorphism of complexes of Zariski sheaves on $\Sm_k$
\begin{equation}\label{eq:CorUDeltaVconginteginvfunct}
\mathrm{Cor}( - \times_k \Delta^\bullet_k , V )) \cong 
\mathbb Z\oplus k[-\times_k D]^\times,
\end{equation}
where
the right side denotes the complex concentrated in degree zero.
\end{theorem}
\begin{proof}
    The proof is similar to \cite[Theorem 4.1]{MVW} with the closed subscheme $D$ instead of the closed subscheme $\{0\}$ in $\PP^1_k$.
\end{proof}

Recall isomorphisms
\begin{equation}\label{eq:HomUYlncongHnisCorDeltaUYGm}
\Hom_{\mathbf{DM}(k)}(U,Y(l)[n])
\cong \mathbb H^{-n+l}_\nis(\mathrm{Cor}(U\times_k\Delta^\bullet_k,Y\otimes \Gm^{\wedge l}))\\
\cong \mathbb H^{-n+l}_\zar(\mathrm{Cor}(U\times_k\Delta^\bullet_k,Y\otimes \Gm^{\wedge l}))
\end{equation}
provided by results of \cite{Voe-hty-inv,Voe-motives,arbitraryclassicalfieldsstr}.

\begin{proposition}\label{prop:extensioncircle}

For any  $l\in \mathbb Z$, there is the isomorphism of presheaves on $\Sm_k$
\begin{equation}\label{eq:prop:anyzZarlochA1A1-z}h^l_\mathrm{mot}(\widetilde T^z)\simeq z_*h^l_\mathrm{mot}(z^*(T)).\end{equation}
\end{proposition}
\begin{proof}
The isomorphism \eqref{eq:prop:anyzZarlochA1A1-z} 
is provided by 
the following quasi-isomorphisms
of complexes of presheaves on $\Sm_k$ 
\begin{multline}\label{eqLLnisCPrA1k-zA1kz-0A1kz-0}
L_\nis\mathrm{Cor}_k(-\times_k\Delta^\bullet_k,\A^1_k-z)
{\simeq}\\
L_\nis z_*\mathrm{Cor}_{k(z)}(-\times_k\Delta^\bullet_{k(z)},\A^1_{k(z)}-0)
{\simeq}\\
z_* L_\nis\mathrm{Cor}_{k(z)}(-\times_k\Delta^\bullet_{k(z)},\A^1_{k(z)}-0),
\end{multline}
where the left quasi-isomorphism follows from \Cref{lm:LA1htr(V)},
and the right one holds because
the functor $z\colon \Sm_{k}\to\Sm_{k(z)}$ preserves 
semi-local henselian essentially smooth schemes.
\end{proof}


\begin{proposition}\label{cor:extensioncircle}
For any $C\in \Sm_k$ such that $\operatorname{dim}_k C=1$,
and any closed reduced zero-dimensional subscheme $D$
,
for any $U\in\Sm_k$,
there is the isomorphism 
\begin{equation}\label{eq:cor:pnszZarlochA1A1-z}
\Hom_{\DM(k)}(U,C/(C-D)[l])\cong
\begin{cases}0, &l\neq -1,0,\\ \mathcal O(U\times D)^\times, &l=-1\\ \mathrm{Pic}(U\times D), &l=0,\end{cases}\end{equation}
see Notation \ref{notation:hmotlsimeqHomDMk} for the left side.
\end{proposition}
\begin{proof}
Because of \Cref{eq:lm:CVveeA1kfz}
the claim reduces to $C=\A^1_k$ 
and $D=z$ 
for some closed point $z$ in $\A^1_k$.
Because of \eqref{eq:prop:anyzZarlochA1A1-z}
the claim reduces to $C=\A^1_k$ and $D=\{0\}$. 
Since by \eqref{eq:HomUYlncongHnisCorDeltaUYGm} 
\[
\Hom_{\DM(k)}(-,\A^1_k/(\A^1_k-0)[l])\cong
\mathbb{H}^{l}_\zar(-,\bbZ(1)[2])
\cong H^{p,1}_\mathrm{mot}(-,\bbZ)\]
for $l=p-2$, see 
\cite[Definitions 3.1 and 3.4]{MVW},
the isomorpihsms \eqref{eq:cor:pnszZarlochA1A1-z} 
follow 
by \cite[Corollary 4.2]{MVW}.
\end{proof}
\begin{remark}\label{rem:imagesofdivisorsofpolynomialsuitleadterm}
    Under the assumptions of \Cref{cor:extensioncircle}
    for a given $U\in\Sm_k$,
    consider a function $f\in k[U\times_k\A^1_k]\cong k[U][t]$ 
    such that 
    $f=t^n+c_{n-1}t^{n-1}+\dots+c_0$
    and $f$ is invertible on $U\times_k D$.
    By the construction of the quasi-isomorphism \eqref{eq:CorUDeltaVconginteginvfunct} it follows
    that
    the morphism in 
    $\Hom(U\times_k\Delta^\bullet_k,V/\pt_k)
    \cong
    \Hom(U\times_k\Delta^\bullet_k,\A^1_k/V[-1])$
    defined by the element $\operatorname{div}(f)$
    corresponds to $f\big|_{U\times D}\in k[U\times D]^\times$.
%
%
%
\end{remark}
\subsection{Homomorphisms from a pair to a pair}

\begin{definition}\label{def:ratinalpoint}
For any scheme $V$ of the form $\A^1_k-z$, 
we define the rational point $p:=\{1\}$ if $\{1\}\in V$, or $p:=\{0\}$ otherwise.
Then the morphism $\pt_k\to V$ given by $p$ defines the decomposition
$V\simeq V/\pt_k\oplus \pt_k$ in $\DM(k)$.
\end{definition}
\begin{remark}\label{rem:ratpointHomoplusdecomposition}
For any $V_0 = \A^1_k-z_0$, $V_1 = \A^1_k-z_1$,
\Cref{def:ratinalpoint} defines
morphisms
\begin{equation}\label{eq:Hom:pointed-unpointed-pointed}
\Hom(V_0/\operatorname{pt}_k,V_1/\operatorname{pt}_k)\to
\Hom(V_0,V_1)\to
\Hom(V_0/\operatorname{pt}_k,V_1/\operatorname{pt}_k).
\end{equation}
Moreover, both morphisms in \eqref{eq:Hom:pointed-unpointed-pointed} preserve the composition in $\DM(k)$
because the composite morphism $\pt_k\xrightarrow{p} V_0\to\pt_k$ is identity, and $\Hom(V_0/\pt_k,V_1)=0$.
\end{remark}

\begin{definition}
Denote by $F$ the set of
regular functions 
\[f\in k[\A^1_k\times_k\A^1_k]=k[\A^1_k][t]\] 
such that 
$f=t^n+c_{n-1}t^{n-1}\cdots+c_0$, where $c_i\in k[\A^1_k]$, and 
$f$ is invertible on $V_0\times_k z_1$.    
\end{definition}

\begin{lemma}\label{sblm:opensubschA1k}
Let 
$V_0 = \A^1_k-z_0$,
$V_1 = \A^1_k-z_1$.
Then 
\begin{equation}\label{eq:HomV0ptV2ptltrivial}
    \Hom(V_0/\pt_k,(V_1/\pt_k)[l])=0
\end{equation}
and there is an isomorphism
\begin{equation}\label{eq:HomVzVocongkVzzokAzo}
\Hom(V_0/\pt_k,V_1/\pt_k)\cong k[V_0\times_k z_1]^\times/k(z_1)^\times
\end{equation}
such that 
\begin{equation}\label{eq:divfmapstofrestr}
\operatorname{div}(f\big|_{V_0\times_k V_1})\mapsto f\big|_{V_0\times_k z_1},\quad\forall f\in F,
\end{equation}
see Notation \ref{notation:divelementinCor} and \ref{notation:cCortoHom} and \eqref{eq:Hom:pointed-unpointed-pointed}. 
\end{lemma}
\begin{proof}
Let $U$ stands for $V_0$ or $\pt_k$.
Then $\Pic(U\times z_1)\cong 0$, 
and it follows
by \Cref{cor:extensioncircle}
that
\begin{equation*}
\Hom(U,V_1/\pt_k) \simeq 
k[U\times_k z_1]^\times
,\quad
\Hom(U,(V_1/\pt_k)[l])\cong 0, \quad l\neq 0,
\end{equation*}
because
$V_1/\pt_k\simeq \A^1_k/V_1[-1]$.
Then \eqref{eq:HomV0ptV2ptltrivial} and \eqref{eq:HomVzVocongkVzzokAzo} follow.
The mapping \eqref{eq:divfmapstofrestr} is provided
by \Cref{rem:imagesofdivisorsofpolynomialsuitleadterm}.
\end{proof}

\begin{lemma}\label{lm:generaorsCorAlOpairs}
Let 
$V_0 = \A^1_k-z_0$,
$V_1 = \A^1_k-z_1$.
The group 
$\Hom(V_0/\pt_k,V_1/\pt_k)$ 
is generated by the elements of the form $\operatorname{div}(f\big|_{V_0\times_k V_1})$ for all $f\in F$, see Notation \ref{notation:divelementinCor} and \eqref{eq:Hom:pointed-unpointed-pointed}. 
\end{lemma}
\begin{proof}
    Let $c\in\Hom(V_0/\pt_k,V_1/\pt_k)$, and let 
    $v\in k[V_0\times_k z_1]^\times$ be an element that class in $k[V_0\times_k z_1]^\times/k(z_1)^\times$ equlas the image of $c$ with respect to \eqref{eq:HomVzVocongkVzzokAzo}.
    Then
    $v=r_0/r_\infty$, for some functions $r_0,r_\infty\in k[\A^1_k\times_k z_1]$
    invertible on $V_0\times_k z_1$.
    For each $i=0,\infty$, 
    for a large enough $n_i\in\bbZ$, there is 
    $f_i\in k[\A^1_k][t]$ 
    such that 
    $f_i=t^{n_i}+c_{n_i-1}t^{n_i-1}\cdots+c_0$, and 
    $f_i$ equals $r_i$ on $\A^1_k\times_k z_1$.
    Then $c$ equals $\operatorname{div}(f_0)-\operatorname{div}(f_\infty)$ because of \eqref{eq:HomVzVocongkVzzokAzo}.
\end{proof}

\begin{proposition}\label{lm:opensubschA1k}
(a)
Let 
$V_0 = \A^1_k-z_0$,
$V_1 = \A^1_k-z_1$.
Consider 
the unique homomorphism
\begin{equation}\label{eq:HomV/pr=Cor}
\Hom(V_0/\pt_k,V_1/\pt_k)\to
\mathrm{Cor}_k(z_0,z_1);\quad
    \operatorname{div}(f\big|_{V_0\times_k V_1})\mapsto \operatorname{div}(f\big|_{z_0\times_k z_1}),\quad\forall
    f\in F,
\end{equation}
see Notation \ref{notation:divelementinCor} and \eqref{eq:Hom:pointed-unpointed-pointed}, 
provided by \Cref{lm:generaorsCorAlOpairs}.
Then \eqref{eq:HomV/pr=Cor} is an isomorphism.

(b)
The composite morphism of
\begin{equation}\label{eq:CorAV0AV1HomV0V1Corz0z1}\Cor_k((\A^1_k,V_0),(\A^1_k,V_1))\to\Hom(V_0/\pt_k,V_1/\pt_k)\xrightarrow{\cong}
\mathrm{Cor}_k(z_0,z_1),
\end{equation} 
see Notation \ref{notation:Corpair} and \ref{notation:CorCorpair(P0P1)} 
and \eqref{eq:Hom:pointed-unpointed-pointed}, 
equals the composite morphism of
\begin{multline}\label{eq:CorkA1kA1ktoCorz0z1}
\Cor_k((\A^1_k,V_0),(\A^1_k,V_1))\to\mathrm{Cor}_k(\A^1_k,\A^1_k)\to\\
\mathrm{Cor}_k(z_0,\A^1_k)\cong\mathrm{Cor}_k(z_0,V_1)\oplus\mathrm{Cor}_k(z_0,z_1)\to\mathrm{Cor}_k(z_0,z_1).
\end{multline}

(c)
    The map \eqref{eq:CorkA1kA1ktoCorz0z1} preserves the composition.
\end{proposition}
\begin{proof}
(a,b)
By \Cref{sblm:opensubschA1k}
$\Hom(V_0/\pt_k,V_1/\pt_k)\cong k[V_0\times_k z_1]^\times/k[\A^1_k\times_k z_1]^\times$.
There is the commutative diagram
\begin{equation}\label{eq:invcycA1open}\xymatrix{
k[V_0\times_k z_1]^\times/k[\A^1_k\times_k z_1]^\times \ar@{^(->}[d] \ar[r]^<<<<<<{\simeq} & Z_0(z_0\times_k z_1)\ar@{^(->}[d]\\
k(\A^1_k\times_k z_1)^\times/k[\A^1_k\times_k z_1]^\times \ar@{->>}[d] \ar[r]^<<<<<<{\simeq} & Z_0(\A^1_k\times_k z_1)\ar@{->>}[d]\\
k(V_0\times_k z_1)^\times/k[V_0\times_k z_1]^\times \ar[r]^<<<<<<{\simeq} & Z_0(V_0\times_k z_1)
.}\end{equation}
The horizontal isomorphisms \eqref{eq:invcycA1open} are induced by the mapping 
\begin{equation}\label{eq:invkAzotoCycAzo}
k(\A^1_k\times_k z_1)^\times \to Z_0(\A^1_k\times_k z_1); \quad v\mapsto \operatorname{div}(v).
\end{equation}
The upper vertical arrows in \eqref{eq:invcycA1open} are injective, and the upper square is pullback.
The right side bottom vertical surjection in \eqref{eq:invcycA1open} is induced by the open immersion $U_0\times_k z_1\to \A^1_k\times_k z_1$, 
and the left side bottom vertical arrow is surjective, because
$k(\A^1_k\times_k z_1)^\times\cong k(V_0\times_k z_1)^\times.$ 
So the isomorphism \eqref{eq:HomV/pr=Cor} follows, because 
\[
\mathrm{Z}_0(z_0\times_k z_1)\cong \mathrm{Cor}_k(z_0,z_1).
\]
We obtained an isomorphism
$\Hom(V_0/\pt_k,V_1/\pt_k)\cong\Cor_k(z_0,z_1)$,
and by the construction
it follows that it equals \eqref{eq:HomV/pr=Cor}.
Finally, 
since
$\Cor_k((\A^1_k,V_0),(\A^1_k,V_1))$
is the subgroup 
in $\Cor_k(\A^1_k,\A^1_k)$
generated by the elements $\operatorname{div}(f)$
for $f\in F$,
the mapping in \eqref{eq:HomV/pr=Cor}
implies that the composite morphisms of \eqref{eq:CorAV0AV1HomV0V1Corz0z1} and \eqref{eq:CorkA1kA1ktoCorz0z1} are equal.

    (c)
    Given $x\in X_{(0)}$, denote by $i_x\colon x\to X$ the closed immersion,
    see Notation \ref{notation:Xpdp}.
    There is the functor
    \begin{equation}\label{eq:XXzmapstozhookrightarrowcoprodXdimkX}
    {\Cor^\mathrm{pair}_k}\to\Cor(\Sch_k)^\mathrm{c.o.pair};\;
    (X,U)\mapsto (X_{(0)},U_{(0)}),
    \end{equation}
    see Notation \ref{notation:Corpair} and \ref{notation:Xpdp},
    because 
    for each 
    $c\in\Cor_k((X,U),(X^\prime,U^\prime))$ and 
    $x$ in $X_{(0)}$, or $U_{(0)}$, 
    there is a finite set of points 
    $x_\alpha^\prime$ in $X^\prime_{(0)}$, or $U^\prime_{(0)}$,
    such that the morphism $c\circ i_x$ in $\Cor(\Sch_k)$ 
    passes throw $\coprod_{\alpha}x_\alpha^\prime$.
    Then the map \eqref{eq:CorkA1kA1ktoCorz0z1} equals the map
    defined by the composite functor of \eqref{eq:XXzmapstozhookrightarrowcoprodXdimkX} and 
    the functor
    \begin{equation}\label{eq:XhookrightarrowXamalghatXmapstohatX}
         \Cor(\Sch_k)^\mathrm{c.o.pair}\to \Cor_k(\Sch_k);\; 
         (X,U)\mapsto X\setminus U.
    \end{equation}
    Hence \eqref{eq:CorkA1kA1ktoCorz0z1} preserves the composition.
\end{proof}

\begin{corollary}\label{cor:prop:CatC/E}
Let $z\in \A^1_k$ be a point such that $k(z)$
is
a simple purely inseparable extension of $k$.
Then the 
composition homomorphism
in the category $\DM(k,\mathbb Z/p\mathbb Z)$, where $p=\operatorname{char}k$,
\[
\begin{array}{ccccccr}
\Hom(T, \widetilde T^z)&\otimes&\Hom(\widetilde T^z,T)&\to& \Hom(T,T)&\cong&\mathbb Z/p\mathbb Z;\\
c_0&\otimes& c_z &\mapsto& c_z\circ c_0
\end{array}\]
is trivial.
%
\end{corollary}
\begin{proof}
By 
\Cref{lm:opensubschA1k} 
since $\operatorname{deg}_k k(z)=0\in \mathbb Z/p\mathbb Z$,
it follows 
the equality 
$c_z\circ c_0=0$ in 
$\Hom(\A^1_k/(\A^1_k-0), \A^1_k/(\A^1_k-0))\cong\mathbb Z/p\mathbb Z$.
\end{proof}

\section{Category of pairs with simple residue field extensions}\label{sect:onedimeclassANDGSphDMB}

Recall that $\Fieldsop^1_k\subset \mathrm{Sch}_k$ is the subcategory spanned by $\operatorname{Spec}(k(\alpha))$
for all simple field extensions of $k$.
Then
the coproduct closure
$(\Fieldsop^1_k)^\amalg$
is the category
of irreducible zero-dimensional finite $k$-schemes that residue fields are simple extensions of $k$.

\begin{proposition}\label{prop:CatC/El}
There is the equivalence of categories
\[\begin{array}{lcl}
(C/V[l]\; |\; (C,V)\in\mathrm{Sm}^\mathrm{pair,1}_k, l\in\bbZ)_{\DM(k)} 
&\simeq& \mathbb Z\times \mathrm{Cor}((\Fieldsop^1_k)^\amalg); \\ 
C/V[l]&\mapsto& (l,(C\setminus V)_\mathrm{red}),
\end{array}\]
where
$\mathrm{Sm}^\mathrm{pair,1}_k$
is the category 
of pairs $(C,V)$, where
$C\in \Sm_k$, $\operatorname{dim}_k C=1$, and $V$ is a dense open subscheme.

\end{proposition}
%
%
\begin{proof}
%
By \Cref{lm:generaorsCorAlOpairs} it follows that 
$\Hom(\A^1_k-z_0,\A^1_k-z_1)$ is generated by 
the morphisms defined by the elements 
in $\Cor_k((\A^1_k,\A^1_k-z_0),(\A^1_k,\A^1_k-z_1))$
for any closed points $z_0$ and $z_1$.
The homomorphism \eqref{eq:CorkA1kA1ktoCorz0z1} 
and the left homomorphism in \eqref{eq:CorAV0AV1HomV0V1Corz0z1}
both respect the composition by 
\Cref{lm:opensubschA1k}(c)
and \Cref{rem:ratpointHomoplusdecomposition}. 
Hence the isomorphism \eqref{eq:HomV/pr=Cor} defines the fully faithful functor
\[
(\A^1_k/(\A^1_k-z)| z\in \A^1_k)_{\DM(k)}
\to \Cor(\Sch_k)
.\]
Then by \Cref{lm:opensubschA1k}(a,b) and \Cref{eq:lm:CVveeA1kfz} the claim follows.
%
%
%
%
\end{proof}

\begin{theorem}\label{th:onedim:subcatDMCorB}

Let $\DM(B)$ be as in Notation \ref{notation:categoriesSH(B)DM(B)}.
For each $d\in \mathbb Z_{\geq 0}$,
there is 
an
equivalence of 
categories 
\begin{equation}\label{eq:simpleqext:subcatDMB(Xx)simeqsubcatCorB(x)}
\begin{array}{lcl}
(\;
X/(X-x)[l]\;|\;
    (X,x)\in\Sm^{d,\cdot,1}_B, l\in \mathbb Z
\;)_{\DM(B)}
&\simeq&
\mathbb Z\times \Cor_B(\Fieldsop^1_B),\\
X/(X-x)[l]&\mapsto& (l,x/B)
,\end{array}
\end{equation}
where 
$x/B\in\Fieldsop_B$ stands for the object given by the morphism $x\to B$, 
see 
Notation \ref{notation:Smddot}, \ref{notation:fsubcat} for the left side,
and
Notation \ref{notation:PointsB},
\ref{notation:catZ} for the right side.
\end{theorem}
\begin{proof}
    By \Cref{prop:CatC/El} the claim holds when $B=\Spec k$ for a field $k$.
    For any base scheme $B$, 
    the claim follows 
    from the result over each 
    closed point $z\in B$
    because of
    the reflective adjunction $z^*\dashv z_*$,
    $z_*\colon \DM(z) \leftrightarrows \DM(B)\colon z^*$
    provided by \cite{Cisinski-Deglise-Triangmixedmotives},
    and 
    since
    for any $X\in \Sm_B$, and $x\in X$ be such that 
    $p(x)\neq z$.
    there is the isomorphism
    $
    z^*(X_x/(X_x-x))\simeq *, 
    $
    of pro-objects in $\mathbf{H}^\bullet(z)$.
\end{proof}

\section{Classification results}
\begin{theorem}\label{th-postintro:classif}

Let
semi-local essentially smooth $B$-schemes 
$X$ and $X^\prime$
with
the subschemes of closed points
$\mathrm{x}=\coprod_{i=0}^n x_i$, 
$\mathrm{x}^\prime=\coprod_{i=0}^{n^\prime} x^\prime_i$
%
be such that
$\dim_B X=\dim_B X^\prime$,
and
for some field extension $K/k$
there is 
an isomorphism
\[K[\mathrm{x}]\cong \prod_{i=0}^{N} K(\alpha_i),\]
where $K(\alpha_i)$ are simple extensions of $K$,
and similarly for $K[\mathrm{x}^\prime]$. 

Then
an isomorphism 
$X/(X-\mathrm{x})\simeq X^\prime/(X^\prime-\mathrm{x}^\prime)$
in 
$\mathbf{SH}(B)$ 
or $\mathbf{DM}(B)$
implies
an isomorphism
$\mathrm{x}\cong \mathrm{x}^\prime\in \Sch_B.$
\end{theorem}\begin{proof}
    The claim follows by \Cref{th:onedim:subcatDMCorB}.
\end{proof}

\begin{lemma}\label{lm-postintro:classif:cls}
Let $X,X^\prime\in\Sm_k$.
Let 
$\mathrm{x}=\coprod_{i=0}^n x_i$
$\mathrm{x}^\prime=\coprod_{i^\prime=0}^{n^\prime} x^\prime_{i^\prime}$
be finite subsets of closed points
such that 
for each $x\in\mathrm{x}$ and $x^\prime\in\mathrm{x}^\prime$,
$\operatorname{dim}^x_B X=d=\operatorname{dim}^{x^\prime}_B X^\prime$
for some $d\in\mathbb Z$.
Suppose there is an isomorphism
\begin{equation}\label{eq:XXxcongTP}
X/(X-\mathrm{x})\cong X^\prime/(X^\prime-\mathrm{x}^\prime)\in \DM(k),
\end{equation}
that induces morphisms
\[c_{i,i^\prime}\colon X/(X-x_i)\to X^\prime/(X^\prime-x^\prime_{i^\prime}).\]
Suppose 
$x^\prime_0\cong\pt_k$,
and
$c_{i,0}=0$
for all 
$i< j$ for some $j\in\{0,\dots,n\}$.

Then 
there is an isomorphism
\begin{equation}\label{eq:xcongptE}
x_{i}\cong \pt_k
\end{equation}
for some 
$i\geq j$.
\end{lemma}

\begin{proof}

    We prove the claim by induction on 
    $v=\sum_{i\in\{j,\dots,n\}} \operatorname{deg}_k k(x_i)$. 
    If $v=0$, 
    the claim is trivial. 
    Suppose the claim is proven for all input data such that $v$ is less then for the give one.
    
    Suppose there is no isomorphism \eqref{eq:xpcongtE}.
    Up to a premutation
    we may assume that
    \begin{equation}\label{eq:minimalityofdegkxzero}
        \operatorname{deg}_k k(x_j)\leq \operatorname{deg}_k k(x_i)
    \end{equation}
    for all $i\in\{j,\dots,n\}$.
    Since $k(x_j)$ is finite over $k$, and $k(x_j)\not\cong k$,
    there is $K\subset k(x_j)$ such that $k(x_j)\cong K(\alpha)$ for some $\alpha\in k(x_j)$.
    Since
    by \Cref{eq:minimalityofdegkxzero} $K$ 
    does not contain any one of $k(x_i)$, $i\in\{j,\dots,n\}$.
    using the base change along the extension $K/k$,
    we reduce the claim to the situation, when $k(x_j)\cong k(\alpha)$.
     
    By assumption
    $X^\prime/(X^\prime-x_0^\prime)\cong T^{\wedge d}$.
    By \Cref{th:onedim:subcatDMCorB} 
    \begin{equation}\label{eq:Homptkx0x0ptk}\begin{array}{lclcl}
         \Hom(T^{\wedge d},X/(X-x_j))&\cong&\mathrm{Cor}(\pt_k,x_j), \\
         \Hom(X/(X-x_j),T^{\wedge d})&\cong&\mathrm{Cor}(x_0,\pt_k).
    \end{array}
    \end{equation}
    Then $c_{0,j}$ and $c_{j,0}$
    are given by elements
    $
    c_{0,j}\in\mathrm{Cor}(\pt_k,x_0)
    $,
    $
    c_{j,0}\in\mathrm{Cor}(\pt_k,x_0)
    $.
    Consider the base change along the extension 
    $k^\mathrm{alg}/k$.
    Consider the image of
    the isomorphism \eqref{eq:XXxcongTP}
    along the equivalence of categories 
    \eqref{eq:simpleqext:subcatDMB(Xx)simeqsubcatCorB(x)} from
    \Cref{th:onedim:subcatDMCorB} applied over $k^\mathrm{alg}$.
    Since $k(x_j)\not\cong\pt_k$ by the assumption,
    we 
    conclude 
    that 
    the base changes of $c_{0,j}$ and $c_{j,0}$
    are trivial. 
    Consequently,
    \[c_{0,j}=0,\quad c_{j,0}=0.\]
    Then  
    by the inductive assumption 
    it follows that 
    $x_i\cong\pt_k$
    for some $i>j$,
    that controdicts to \eqref{eq:minimalityofdegkxzero}.
\end{proof}

\begin{theorem}\label{th-postintro:classif:cls}

Let $X,X^\prime\in\Sm_B$,
and let $x\in X$, $x^\prime\in X^\prime$ be closed points.
The following conditions are equivalent:
\begin{itemize}
    \item[(1)] $\mathrm{x}\cong \mathrm{x}^\prime\in \Sch_B$,
    and
    $\dim^{x}_{B}X=\dim^{x^\prime}_B X^\prime$
    for each $x\in X$,
    and
    corresponding $x^\prime\in X^\prime$;
    \item[(2)] $X/(X-{x})\simeq X^\prime/(X^\prime-{x}^\prime)$
    in
    $\mathbf{H}^\bullet(B)$;
    \item[(2')] $X/(X-{x})\simeq X^\prime/(X^\prime-{x}^\prime)$
    in
    $\mathbf{DM}(B)$;
\end{itemize}
\end{theorem}
\begin{proof}
    (1) implies (2) by \Cref{fact:ClassifParta} and \Cref{rem:Nisnevicheq}.
    (2) implies (2').
    We are going to prove that (2') implies (1).
    Similarly as in \Cref{th:onedim:subcatDMCorB}
    using base changes along morphisms 
    $z\to B$ for all points $z\in B$
    we reduce the claim 
    to the base field setting.
    Consider the base change to the separable 
    closure of the base field $k$.
    Then equivalence $X/(X-{x})\simeq X^\prime/(X^\prime-{x}^\prime)$ 
    induces an equivalence
    \begin{equation}\label{eq:sumsofTwsomed}
        \bigoplus_{\operatorname{sdeg}_k k(x)} T^{\wedge \dim^x_k X}
        \cong 
        \bigoplus_{\operatorname{sdeg}_k k(x^\prime)} T^{\wedge \dim^{x^\prime}_k X}.
    \end{equation}
    Hence $\dim^x_k X=\dim^{x^\prime}_k X$.

    Without loss of generality we may assume that 
    \begin{equation}\label{eq:th-postintro:classif:cls:degkx0lgeqdegxii}
        \operatorname{deg}_k k(x)\leq \operatorname{deg}_k k(x^\prime).
    \end{equation}
    Consider the base change 
    with respect to the field extension $K/k$, $K=k(x)$. 
    Applying \Cref{lm-postintro:classif:cls} over $K$,
    we conclude that there is 
    a morphism 
    $x\to x^\prime\times_k \Spec K$
    in $\Sch_{k}$,
    and conequently,
    a morphism 
    $x\to x^\prime$
    in $\Sch_{k}$.
    Then by \eqref{eq:th-postintro:classif:cls:degkx0lgeqdegxii},
    it follows that there is an isomorphism
    \[x\cong x^\prime.\]
\end{proof}

\appendix

\section{Cycle modules} 
\label{sect:CycleModules}
In \cite{Rost1996}, M. Rost
defined the category of cycle modules $\mathcal{MC}ycle(k)$
over a field $k$.
It is proven
in \cite{DegliseModuleshomotopiques}
by F. D{\'e}glise 
that 
for a perfect fields $k$ or after the inverting of the exponential characteristic of $k$ in $\Lambda$
there is
an equivalence
    \begin{equation}\label{eq:MmodDMheart}
    \mathcal{MC}ycle(k,\Lambda)
    \simeq
    \mathbf{DM}(k,\Lambda)^{\heartsuit}
    \end{equation}
    where $\heartsuit$ relates to the homotopy t-structure heart.
\begin{theorem}\label{eq:hearttshomottructure}
    Equivalence \eqref{eq:MmodDMheart}
    holds for any field $k$ with integral coefficients.
\end{theorem}
\begin{proof}
    Let us remind 
    the Rost transform 
    functor $\mathcal{R}\colon \mathcal{MC}ycle(k)\to\mathbf{DM}(k)$
    from \cite{DegliseModuleshomotopiques}.
    For a given
    $F\in \mathcal{MC}ycle(k)$,
    there is a sheaf of abelian groups on $\Sm_k$ such that 
    \begin{equation}\label{eq:R(F)(U)}
    \mathcal{R}(F)(U)
    \cong 
    \operatorname{Ker}
    (F(U^{(0)})\to \bigoplus_{z\in U^{(1)}}F(z))
    \end{equation}
    for
    an essentially smooth local henselian $U$,
    that is $\A^1$-invariant by the axioms of cycle modules  
    and has $\mathrm{Cor}(k)$-transfers by \cite{zbMATH05000059}.
    %
    By \cite[Th 4.1.6]{déglise2022perversehomotopyheartmwmodules}
    there is the adjuntion
    $\mathcal{R}\colon 
        \mathcal{MC}ycle(k)
    \simeq
    \mathbf{DM}(k)
    \colon \mathbf{H}$
    and $\mathcal R\circ \mathbf{H}\simeq\mathrm{Id}$,
    where
    for a given $C\in\mathbf{DM}(k)$,
    $
    \mathbf{H}(C)(x,c)
    \cong 
    H^{c}_{x}(X,C),
    $
    where
    $c=\codim_{X}x$.
    Since 
    $E=\mathcal{R}(F)$ is $\A^1$-invariant,
    it is strict $\A^1$-invariant by
    \cite{StrictA1invariancefields},
    and by
    \cite[Theorem 4.37]{Voe-hty-inv}
    the Cousin complex 
    \[E(U^{(0)})\to\dots\to\bigoplus_{z\in U^{(c)}}E^c_z(U)\to\dots\]
    is acyclic, where 
    $E^c_z(U)=H^c_z(U_z,E).$
    Then 
    $
    \mathbf{H}\circ\mathcal{R}\simeq\mathrm{Id}$.
\end{proof}

\section{Field generators and scheme dimension}
\label{sect:fgANDsd}

\begin{lemma}\label{lm:Xprime}
Let $p\colon X\to B$ be a smooth morphism of schemes, and $x\in X$.
Denote $z=p(x)\in B$, $k = \mathcal O_{z}(z)$, and $K = \mathcal O_x(x)$. 
Suppose that $K/k$ is finite, and $K\cong k(\phi_1,\dots,\phi_n)$.
Then there are a Zariski neighbourhood $U$ of $x$ in $X$ and a smooth closed subscheme $X^\prime$ in $U$ 
such that $x\in X^\prime$ and $\dim^x_B X^\prime = n$.
\end{lemma}
\begin{proof}
    Since $x$ has an affine Zariski neighbourhood in $X$,
    without loss of generality we can assume that $X$ is affine.
    Since the extension $K/k$ is finite, 
    $x$ is a closed point in $X\times_B z$.
    Denote $d=\dim^x_B X=\dim^x_k{(X\times_B z)}$.
    Let 
    \[\{f_1,\dots,f_{d}\in \mathcal O_X(X)\}\]
    be a set of functions on $X$ 
    such that 
    $f_{i}\big|_{x}=\varphi_i$ for $i\leq n$,
    $f_{i}\big|_{x}=0$ for $i> n$,
    and 
    the differentials $\{d^x_{X\times_B z} f_1,\dots,d^x_{X\times_B z} f_d\}$ at the point $x$ on $X\times_B z$
    are linearly independent.
    Then
%
    the morphism 
    \[f=(f_1,\dots,f_{d})\colon X\to\A^d_B\]
    induces the surjection of $K$-vector spaces 
    \[T_{X\times_B z,x}\to T_{\A^d_k,v},\quad
    v=f(x).\]
    Hence $f$ is \'etale over $x$.
    Let $U$ be an open neighbourhood of $x$ in $X$ such that $f\big|_U$ is \'etale,
    and
    put \[X^\prime := Z(f_{n+1}\big|_U,\dots,f_{d}\big|_U)\subset U.\]
    Then $x\in X^\prime$, $X^\prime\in\Sm_B$, $\dim^x_B X^\prime=n$.
\end{proof}

\begin{lemma}\label{lm:Xprime:XX-x_simeq_XpXp-x_wedge_Td-1}
Let $X\in \Sm_B$, $x\in X$.
For any smooth closed subscheme $X^\prime$ in $X$ such that $x\in X^\prime$,
there is an isomorphism in $\mathbf{H}^\bullet(B)$
\[X_x/(X_x-x)\simeq X^\prime_x/(X^\prime_x-x)\wedge T^{\dim^x_B X-\dim^x_B X^\prime}.\]
\end{lemma}
\begin{proof}
    The claim follows because the normal bundle $N_{X^\prime/X}$ is locally trivial. 
\end{proof}

\section{Presentation of one-dimensional pairs}\label{sect:GabLmconseq}

We use the following consequences of 
Gabber's Presentation Lemma provided by 
\cite[Lemma 3.1]{Gab}, \cite[Theorem 3.2.1]{CTHK}, and \cite{HK}.
\begin{lemma}\label{lm:GabbersPL}
Let $C\in\Sm_k$, $\dim_k C=1$.
For any closed point $z$ in $C$, 
there is 
an \'etale morphism $f\colon C\to \A^1_k$ that maps $z$ isomorphically on $f(z)$.
\end{lemma}
\begin{proof}
    The claim follows by
    \cite[Lemma 3.1]{Gab}, \cite[Theorem 3.2.1]{CTHK}, \cite{HK}. 
\end{proof}

\begin{lemma}\label{eq:lm:CVveeA1kfz}
    For any $C\in \Sm_k$, $\dim_k C=1$, and a dense open subscheme $V=C-D$, 
    \[C/V\cong \bigoplus_{z\in D} \A^1_k/(\A^1_k-v_z)\cong \bigoplus_{z\in D} (\A^1_k-v_z)/\pt_k[1],\]
    for some family of closed points $v_z\in \A^1_k$ such that $v_z\cong z$ for all points $z$ in $D$.
\end{lemma}
\begin{proof}
For any $C$ and $V=C-D$ as above, 
$C/V\cong \bigoplus_{z\in D} C/(C-z)$.
For each $z$, by \Cref{lm:GabbersPL} there are isomorphisms
$
C/(C-z)\simeq \A^1_k/(\A^1_k-f(z))\simeq (\A^1_k-f(z))/\pt_k[1].
$
So the claim follows.
\end{proof}

\printbibliography
\end{document}